\documentclass[final,authoryear,5p,twocolumn]{elsarticle}
\usepackage{tabularx}
\usepackage{graphics}
\usepackage{graphicx}
\usepackage{epsfig}
\usepackage{amssymb}
\usepackage{amsthm}
\usepackage{amsmath,amsfonts}
\usepackage{array}
\usepackage{multirow}
\usepackage{algorithm,float}
\usepackage{float}
\usepackage{algorithmic}
\usepackage{mathrsfs}
\usepackage{epsfig,epstopdf}
\usepackage{textcomp}
\usepackage{booktabs}
\usepackage[colorlinks,linkcolor=blue,citecolor=blue]{hyperref}

\usepackage{hyperref}
\usepackage{bm}
\usepackage{bbm}
\usepackage{color}
\usepackage[title]{appendix}
\usepackage{pifont}
\usepackage{makecell}
\usepackage{url}
\usepackage[nooneline,flushleft]{caption2}
\usepackage{xr}
\usepackage{multirow}
\usepackage{indentfirst}
\biboptions{sort&compress}
\newtheorem{theorem}{Theorem}
\newdefinition{definition}{Definition}
\newdefinition{lemma}{Lemma}
\newdefinition{assumption}{Assumption}
\newdefinition{remark}{Remark}
\newtheorem{corollary}{Corollary}

\allowdisplaybreaks[4]

\newcommand{\bsx}{\boldsymbol{x}}

\newcommand{\bsX}{\boldsymbol{X}}

\newcommand{\sumn}{\sum\limits_{i=1}^n}

\newcommand{\sumT}{\sum\limits_{t=1}^T}

\journal{Automatica}
\allowdisplaybreaks[4]
\begin{document}

\begin{frontmatter}


\title{Dynamic Regret of Distributed Online Frank-Wolfe Convex Optimization\tnoteref{footnoteinfo}}
\tnotetext[footnoteinfo]{
This paper was not presented at any conference.}

\author[a1]{Wentao~Zhang}\ead{iswt.zhang@gmail.com}
\author[a2]{Yang~Shi}\ead{yshi@uvic.ca}
\author[a1]{Baoyong~Zhang
\corref{cor}
}
\ead{baoyongzhang@njust.edu.cn}
\cortext[cor]{Corresponding author.}
\author[a1]{Deming~Yuan}\ead{dmyuan1012@gmail.com}
\address[a1]{School of Automation,  Nanjing University of Science and Technology, Nanjing 210094, Jiangsu, P. R. China}
\address[a2]{Department of Mechanical Engineering, University of Victoria, Victoria, BC,  Canada, V8W 3P6}

\begin{abstract}
 This paper considers distributed online convex constrained optimization, in which various agents in a multi-agent system  cooperate to minimize a global cost function through communicating with neighbors over a time-varying network.
 When the constraint set of optimization problem is high-dimensional and complicated, the computational cost of the projection operation  often becomes prohibitive. To handle this problem, we develop a distributed online Frank-Wolfe optimization algorithm combining with gradient tracking technique.  We rigorously establish the dynamic regret bound of the proposed optimization algorithm as $\mathcal{O}(\sqrt{T(1+H_T)}+D_T)$, which explicitly depends on the iteration round $T$, function variation $H_T$, and gradient variation $D_T$.   Finally, the theoretical results are verified and compared in the case of distributed online ridge regression problems.
\end{abstract}


\end{frontmatter}

\section{Introduction}
Recently, distributed optimization over multi-agent network has attracted much attention  due to its wide applications, such as power systems, sensor networks, machine learning, etc (see, e.g., \cite{nedic2018distributed}, \cite{yang2019survey},  \cite{li2022survey}, \cite{liu2020unitary}, \cite{li2021distributed}, \cite{xu2016distributed}, \cite{yuan2022distributedauto}, \cite{9992989}).  In such an optimization problem, distributed online optimization
 can be described as a repeated game as follows:

\begin{enumerate}
  \item At every round $t$, agent $i$ first generates a decision $\bsx_{i,t} \in \boldsymbol{X}$.
   \item Agent $i$ suffers a loss $f_{i,t} (\bsx_{i,t})$ and the adversary reveals  the information about loss function $f_{i,t}$.
  \item Then, agent $i$ uses the information about loss function $f_{i,t}$  to construct the next decision $\bsx_{i,t+1}$.
\end{enumerate}
The main task of the network is  to minimize the sum   of local objective function in  problem (\ref{problem definition}) by mean of information exchange among all agents and local computation, where every agent only knows the information about itself and its neighbors at each round.
\begin{flalign}\label{problem definition}
\min\limits_{\boldsymbol{x}\in \boldsymbol{X}}\,\,\sum\limits_{t=1}^T F_t (\boldsymbol{x})
\end{flalign}
where $F_t (\boldsymbol{x})= \sum_{i=1}^n{f_{i,t}}( \boldsymbol{x} )$, $\boldsymbol{X}$ is convex and compact set in  ${\mathbb{R}}^d$, and the function $f_{i,t}$ is convex in   $\bsX$.  To measure the performance of the proposed algorithm, the dynamic regret  $\textbf{Regret}_{d}^j (T)$ is defined, which represents the total  sum  over time $T$  of the difference between the cumulative cost $F_t { (\boldsymbol{x}_{j,t})}$ of the agent $j$ and the cumulative cost at the optima $\bsx_t^*$.
\begin{flalign}\label{Regret-j}
\textbf{Regret}_{d}^j (T) =   \sum_{t=1}^T F_t { (\boldsymbol{x}_{j,t})} -  \sum\limits_{t=1}^T F_t  (\boldsymbol{x}_t^*)
\end{flalign}
where $\boldsymbol{x}_t^* \in {\arg\min}_ {\boldsymbol{x} \in \boldsymbol{X}}    \ F_t(\boldsymbol{x}).$

 In such distributed online (off-line) optimization problems with constraint sets, projection operations are usually used as a fundamental technique to deal with constraints,
 such as distributed online gradient descent in \cite{sundhar2010distributed}.
 In general, the projection step
 is equivalent to solving a convex quadratic problem \cite{hazan2012projection}. However,  in high-dimensional and complex constrained optimization problems such as multiclass classification in \cite{pmlr-v70-zhang17g}, optimal control in \cite{wu1983conditional}, matrix completion in \cite{hazan2012projection,7883821}, electric vehicle charging in \cite{7496921}, and semidefinite programs in  \cite{hazan2008sparse}, projection operations cause a heavy computational burden.
In contrast, Frank-Wolfe (FW) method avoids such operations with expensive computational cost  through solving a linear minimization oracle.
 \begin{table*}[t]
 \scriptsize
\begin{center}
  \caption{The comparison among relevant works on online FW convex optimization.}
 \begin{tabular}{ccccccc}
 \toprule
 Reference & Loss function&Distributed & \makecell[c]{Dynamic\\ regret }&\makecell[c]{Dynamic\\ network} &\makecell[c]{Linear \\oracle } & Regret bound\\    \midrule
 \cite{pmlr-v70-zhang17g} & Convex and continuous &\ding{51} & \ding{55}  & \ding{55} &  $\mathcal{O}(T)$ & $\mathcal{O}(T^{3/4})$\\
 \cite{pmlr-v119-wan20b} & Convex and continuous &\ding{51} & \ding{55} &\ding{55} &$\mathcal{O}(\sqrt{T})$ & $\mathcal{O}(T^{3/4})$\\
 \cite{wan2021projection} & Strongly convex and continuous&\ding{51}  & \ding{55} &\ding{55} &$\mathcal{O}(T^{1/3})$& $\mathcal{O} (T^{2/3}\log T)$ \\
 \cite{thuang2022stochastic} & Convex and smooth &\ding{51} & \ding{55} &\ding{55} &$\mathcal{O} (T^{3/2})$& $\mathcal{O} (\sqrt{T})$ \\
 \multirow{2}{*}{\cite{kalhan2021dynamic}} & \multirow{2}{*}{Convex and smooth} & \multirow{2}{*}{\ding{55}} &  \multirow{2}{*}{\ding{51} }&  \multirow{2}{*}{/} & $\mathcal{O}(T)$ & $\mathcal{O}\left(\sqrt{T}\left(1+H_T+\sqrt{D_T}\right)\right)$\\
                          &  & & &  & $\mathcal{O}(T^{3/2})$ & $\mathcal{O}\left(1+H_T+\sqrt{T}\right)$\\
\multirow{2}{*}{ \cite{Wan_Xue_Zhang_2021}} & Convex and continuous  &\ding{55} & \ding{51} & / &$\mathcal{O}(T\log_2 T)$ & $\mathcal{O} \left(\max\left\{\sqrt{T}, T^{2/3}H_T^{1/3}\right\}\right)$ \\
 &Strongly convex and continuous  &\ding{55} & \ding{51} & / &$\mathcal{O}(T\log_2 T)$& $\mathcal{O} \left(\max\left\{\sqrt{TH_T\log T}, \log T\right\}\right)$ \\
  This work & Convex and smooth &\ding{51} &  \ding{51} & \ding{51} &$\mathcal{O} (T)$ & $\mathcal{O}\left(\sqrt{T(1+H_T)}+D_T\right)$\\
\bottomrule
\end{tabular}
 \label{table_1}
\end{center}
\vspace{-2.5em}
\end{table*}
 Because of its low computational cost, FW method has been widely utilized  in distributed online (off-line) optimization in recent years.
In \cite{pmlr-v70-zhang17g}, the authors  earlier proposed an online distributed FW algorithm by extending online centralized FW algorithm in \cite{hazan2012projection}  and obtained the static regret bound $\mathcal{O} (T^{3/4})$.  \cite{pmlr-v119-wan20b} considered  an improved variant under full-information feedback and bandit feedback. Based on the idea of dividing time $T$ into $\sqrt{T}$ equally-sized blocks, the frequency of communication between agents was reduced and the  related regret upper bounds $\mathcal{O} (T^{3/4})$ and $\mathcal{\tilde{O}} (T^{3/4})$ were established, respectively. The paper \citep{wan2021projection} further exploited the improved convergence results  under the condition of strong convexity on the basis of \cite{pmlr-v119-wan20b}.  \cite{thuang2022stochastic} analyzed two algorithm versions of exact and stochastic gradient under smooth loss function and showed the regret upper bound $\mathcal{O}(\sqrt{T})$. However,  in both algorithms, the extra step size loop for each agent slows down the computation at time $t$.

\allowdisplaybreaks {Up to now, there is few research work on dynamic regret of online FW algorithms, especially in distributed scenarios. Dynamic regret is a more stringent and effective performance metric than static regret because of its dynamic rather than fixed benchmark, whose bound is generally related with the regularity of the optimization problem.   From \cite{besbes2015non}, it is well known that dynamic regret can not achieve sublinear convergence unless the variation budget satisfies  sublinear in $T$. With that in mind, the \emph{function variation} $H_T$ and  \emph{gradient variation} $D_T$  related to the bound of $\textbf{Regret}_d^j (T)$  are defined as
\begin{flalign}
H_T&= \sumT \max\limits_{i \in \mathcal{V}} \max\limits_{\bsx \in \boldsymbol{X}} |  f_{i,t+1}(\bsx)- f_{i,t}(\bsx)|,  \label {H_T}\\
D_T&= \sumT \max\limits_{i \in \mathcal{V}} \max\limits_{\bsx \in \boldsymbol{X}}\left\| \nabla f_{i,t+1}(\bsx)-  \nabla f_{i,t}(\bsx)\right\| \label {D_T}.
\end{flalign}
In \cite{kalhan2021dynamic}, the authors analyzed the dynamic regret bounds of several centralized online FW algorithms, in which $H_T$ has a limitation about application range. In this paper, we aim to further improve the range under the same conditions, thus the proposed algorithm has stronger applicability.  \cite{Wan_Xue_Zhang_2021} considered a novel centralized online FW algorithm using a restarting strategy. In detail, the comparisons among relevant works on online FW convex optimization are summarized in Table \ref{table_1}.
}

It's not hard to notice that the works in \cite{kalhan2021dynamic} and \cite{Wan_Xue_Zhang_2021} only analyze the dynamic regret  for centralized online optimization.   It is well known that the large-scale optimization problems are difficult to be addressed by the centralized online algorithms due to the computational bottleneck of single machine. From the comparisons and analysis in Table 1, the dynamic regret of online FW algorithm under distributed scenarios needs to be developed. Hence,  the \underline{d}istributed \underline{o}nline \underline{F}rank-\underline{W}olfe \underline{c}onvex \underline{o}ptimization algorithm (DOFW-CO) is designed to fill this gap in this paper.
This paper makes the following contributions.

Firstly, we develop a distributed online Frank-Wolfe convex optimization algorithm that can efficiently deal with the high-dimensional and complicated  constraint set, which alleviates the high computational burden imposed by the projection operator. Meanwhile, the gradient tracking technique is utilized in Algorithm DOFW-CO  to update the gradient change of loss function by using history information.  Moreover, different from the previous communication topology in \cite{pmlr-v70-zhang17g}, \cite{wan2021projection}, \cite{pmlr-v119-wan20b}, \cite{thuang2022stochastic}, the distributed optimization algorithm over a time-varying network topology is developed, which is more practical and general than static network.

Secondly, inspired by \cite{kalhan2021dynamic}, the dynamic regret bound $\mathcal{O}(\sqrt{T(1+H_T)}+D_T)$ for the proposed algorithm is established, which can recover the centralized result and further broaden the range of $H_T$ compared to \cite{kalhan2021dynamic} under the same complexity of linear oracle. The dynamic regret analysis of the FW method in a distributed scenario, for the first time, is developed in this work.  Finally, the case of distributed online ridge regression problems is simulated to verify the performance of the proposed algorithm.

\textbf{Notation}: ${\mathbb{R}}^n$  represents  the Euclidean space with $n$ dimensions.  $\mathbb{Z}$ $(\mathbb{Z}_+)$ represents the (positive) integers set.  The Euclidean  norm  of a vector $\boldsymbol{z}$ is denoted as $\|\boldsymbol{z}\|$. $[A_t]_{ij}$ signifies the element in the $i$-th row and $j$-th column of matrix $A_t$ and $[\boldsymbol{w}]_{i}$ denotes the $i$-th element of vector $\boldsymbol{w}$.

\section{Problem Formulation}
\subsection{The Optimization Problem}

 Let $\mathcal{G}_t=\{\mathcal{V},\mathcal{E}_t,A_t \}$ represent a directed time-varying network with  the set $\mathcal{V} : = \{1,\ldots,n\} $ of agents, the edge set $ \mathcal{E}_t \subseteq \mathcal{V} \times \mathcal{V} $ and the weighted adjacency matrix $A_t \in \mathbb{R}^{n \times n}$. In the network, agent $i$ has permission to communicate with the agents on inner neighbor sets $\mathcal{N}_{i}^{\text {in }}(t)=\{j \mid(j, i) \in \mathcal{E}_t\} \cup\{i\}$ of agent $i$.  Further, $[A_{t}]_{ij}>0$ holds when $j \in \mathcal{N}_{i}^{\text {in }}(t)$, and $[A_{t}]_{ij}=0$ holds otherwise.

The objective of this paper is to design a distributed online algorithm for problem  (\ref{problem definition}) to ensure that the dynamic regret of every agent $j \in \mathcal{V} $  grows sublinearly, i.e., $\lim_{T \rightarrow \infty} (\textbf{Regret}_d^j (T) / T) =0, \forall j \in \mathcal{V}$.
Around the network $\mathcal{G}_t$, the constraint set and loss function  in problem \ref{problem definition}, the following assumptions are made.
\begin{assumption}\label{assump: network}
(a)  There exists a positive scalar $ \zeta $ such that $ [A_{t}]_{ij}>\zeta, t \in \{1,\ldots, T\}$ when $[A_{t}]_{ij}>0$.

(b) $A_t$ satisfies $\sum_{j=1}^{n}  [A_{t}]_{ij}=\sum_{i=1}^{n} [A_{t}]_{ij}=1$ for any $t\in \{1,\ldots, T\}$ and all $i,j \in \mathcal{V }$.

(c)  With some $Q \in \mathbb{Z}_{+}$, the graph's union $\bigcup_{i=k Q+1}^{(k+1)Q} \mathcal{G}_{i}$ is strongly connected for every integer $k \geq 0$.
\end{assumption}


\begin{assumption}\label{assump: ball}
The constraint set $\bsX$ has a finite diameter $M$, i.e., for  $\forall \boldsymbol{x}_1,\boldsymbol{x}_2 \in \boldsymbol{X}, \max_{\boldsymbol{x}_1,\boldsymbol{x}_2 \in \boldsymbol{X}} \| \boldsymbol{x}_1 -\boldsymbol{x}_2 \| \leq M.$
\end{assumption}
\begin{assumption}\label{assump: lips funon}
(Lipschitz Function)
 The function $f_{i,t}$ is $L_{X}$-Lipschitz,  i.e., for $\forall \boldsymbol{x}_1,\boldsymbol{x}_2 \in \boldsymbol{X}$, $
 |f_{i,t}(\boldsymbol{x}_{1})-f_{i,t}(\boldsymbol{x}_{2})| \leq L_{X} \|\boldsymbol{x}_{1}-\boldsymbol{x}_{2}\|, $
where $L_{X}$  is known positive constant.
\end{assumption}
\begin{assumption}\label{assump: lips grad}
(Lipschitz Gradient) The gradient $\nabla f_{i,t} (\boldsymbol{x})$ is $G_{X}$-Lipschitz, i.e., $
\| \nabla f_{i,t} (\boldsymbol{x}_1) - \nabla f_{i,t} (\boldsymbol{x}_2) \|\leq {G_X} \| \boldsymbol{x}_1-\boldsymbol{x}_2\|, \forall \boldsymbol{x}_1,\boldsymbol{x}_2 \in \boldsymbol{X}.$
\end{assumption}
\begin{remark}
In centralized and distributed optimization, Assumptions \ref{assump: network}, \ref{assump: ball}, \ref{assump: lips funon} are standard and similar settings can be seen in \cite{nedic2008distributed, yi2020distributed}, \cite{7883821,besbes2015non}.
According to Lemma 2.6 in \cite{shalev2011online}, Assumption \ref{assump: lips funon} implies  $\|\nabla f_{i,t}(\bsx) \|\leq L_X$.  Assumption \ref{assump: lips grad}  is equivalent to the fact
\begin{flalign}
f_{i,t}(\boldsymbol{x}_{1})-f_{i,t}(\boldsymbol{x}_{2}) & \leq \langle \nabla f_{i,t} (\bsx_2), \bsx_1-\bsx_2\rangle\nonumber \\
 & \quad + \frac{G_X}{2} \|\boldsymbol{x}_{1}-\boldsymbol{x}_{2}\|^2, \forall \bsx_1, \bsx_2 \in \bsX.
\end{flalign}
\end{remark}
\section{Algorithm Design and Convergence Analysis}
\subsection{Algorithm DOFW-CO}
In this section, we first develop Algorithm DOFW-CO. The algorithm description is presented in Algorithm \ref{algorithm 1}. Specifically, the key ingredients include: 1) the gradient tracking and Frank-Wolfe methods are utilized; 2) the gradient $\widehat{\nabla} f_{i,t}$ after the gradient tracking step replaces the traditional gradient $\nabla f_{i,t}  (\hat{\boldsymbol{x}}_{i,t})$ in the linear oracle of Frank-Wolfe step.
\begin{algorithm}[t]
	\renewcommand{\algorithmicrequire}{\textbf{Initialize:} }
	\caption{ (DOFW-CO) Distributed Online Frank-Wolfe Convex Optimization }
	\label{algorithm 1}
	\begin{algorithmic}[1]
		\REQUIRE Initial  variables $\boldsymbol{x}_{i,1} \in \boldsymbol{X}$ and parameter $0 < \alpha \leq 1.$
		
		\FOR {$t=1,2,\cdots,T$}
             \FOR {Each agent $i \in \mathcal{V}$}
             \STATE Agent $i$ receives $\boldsymbol{x}_{j,t}$ from $j \in \mathcal{N}_{i}^{\text {in }}(t)$, and  updates

                \begin{center}
                \quad $\boldsymbol{\hat{x}}_{i,t}= \sum_{j \in \mathcal{N}_{i}^{in}(t)}{[A_t]_{ij} \boldsymbol{x}_{j,t}}.$
                           \setlength{\parskip}{0.4em}
                \end{center}
\STATE  The gradient value $\nabla f_{i,t}  (\boldsymbol{\hat{x}}_{i,t})$  is revealed and agent $i$ executes gradient tracking steps:
            \IF{$t=1 $}
            \STATE $\overline{\nabla} f_{i,1} =  \nabla f_{i,1}  (\hat{\boldsymbol{x}}_{i,1})$,
            \ELSE
            \STATE $\overline{\nabla} f_{i,t} = \widehat{\nabla} f_{i,t-1} + \nabla f_{i,t}  (\hat{\boldsymbol{x}}_{i,t}) - \nabla f_{i,t-1}  (\hat{\boldsymbol{x}}_{i,t-1}).$
            \ENDIF

               \begin{center}
               \quad $\widehat{\nabla} f_{i,t}= \sum_{j \in \mathcal{N}_{i}^{in}(t)}{[A_{t}]_{ij}} \overline{\nabla} f_{j,t},$
               \end{center}

		\STATE Frank-Wolfe step: update

             \begin{center}
               \quad $\boldsymbol{v}_{i,t} = \underset{\boldsymbol{x} \in \boldsymbol{X}}{\arg\min} \left<\boldsymbol{x}, \widehat{\nabla} f_{i,t} \right>, $

               \quad $\boldsymbol{x}_{i,t+1}= \hat{\boldsymbol{x}}_{i,t} +\alpha (\boldsymbol{v}_{i,t}- \hat{\boldsymbol{x}}_{i,t}).$
              \end{center}

		\ENDFOR
		\ENDFOR
	\end{algorithmic}
\end{algorithm}

\subsection{Main Convergence Results}
In this section, the upper bound of dynamic regret defined in (\ref{Regret-j}) for Algorithm \ref{algorithm 1} is analyzed in detail. In order to facilitate the proof and analysis, we define the running average vectors $\bsx_{avg,t}$ and $\boldsymbol{v}_{avg,t}$, the max function variation $f_{t,sup}$ at time $t$, the max gradient variation $g_{t,sup}$ at time $t$ and gradient difference  $\boldsymbol{\delta}_{i,t}$ of agent $i$ as follows:
\begin{gather}\label{average-delta}
\left\{\begin{array}{rcl}
\boldsymbol{x}_{avg,t}&=&\frac{1}{n} \sum\limits_{i=1}^n \boldsymbol{x}_{i,t},
\boldsymbol{v}_{avg,t}=\frac{1}{n} \sum\limits_{i=1}^n \boldsymbol{v}_{i,t} \\
f_{t,sup}&=& \max\limits_{i \in \mathcal{V}} \max\limits_{\bsx \in \boldsymbol{X}} |   f_{i,t+1}(\bsx)-   f_{i,t}(\bsx)| \\
g_{t,sup}&= &\max\limits_{i \in \mathcal{V}} \max\limits_{\bsx \in \boldsymbol{X}}\left\|  \nabla f_{i,t+1}(\bsx)-  \nabla f_{i,t}(\bsx)\right\| \\
\boldsymbol{\delta}_{i,t}&= &\nabla f_{i,t}  (\hat{\boldsymbol{x}}_{i,t}) - \nabla f_{i,t-1}  (\hat{\boldsymbol{x}}_{i,t-1})
\end{array}\right.
\end{gather}
\begin{theorem}\label{theorem 1}
Let the decision sequence $\{ \bsx_{i,t}\}$ be generated by Algorithm \ref{algorithm 1} and suppose Assumptions \ref{assump: network}-\ref{assump: lips grad} hold. Then,     for $T\geq 2$ and $j \in \mathcal{V}$,  the regret is bounded as follows:
\begin{flalign}\label{Eq Theorem1}
 \mathbf{Regret}_d^j(T)
&\leq C_1  +C_2 \alpha T +\frac{2n}{\alpha} H_T + \frac{C_3}{\alpha}  +C_4 D_T
\end{flalign}
where
\begin{align*}
& \ C_1=n L_{X} \sum_{i=1}^n \|\boldsymbol{x}_{i,1}-\boldsymbol{x}_{avg,1}\| + \frac{ 2M n \Gamma}{1-\sigma}\sum_{i=1}^n \|   {\nabla} f_{i,1} (\hat{\boldsymbol{x}}_{i,1}) \|
 \\
 & \quad +\left(nL_X+2MG_X+\frac{4Mn\Gamma G_X}{1-\sigma}\right)\frac{n\Gamma}{1-\sigma} \sum_{i=1}^n  \| {\boldsymbol{x}}_{i,1} \|,  \\
 &C_2=2n^2 L_X M +\left(4M G_X+n L_X +\frac{4Mn\Gamma G_X}{1-\sigma}\right)\frac{ n^2 M  \Gamma}{1-\sigma}\\
 & \quad +\frac{nG_X M^2}{2},
 C_3=n L_X M,  C_4=\frac{2M n^2 \Gamma}{1-\sigma}+n M,\\
   &\sigma=\left(1-{\zeta}/ { 4 n^{2}}\right)^{1 /  Q},\quad  \Gamma=\left(1-{\zeta} //{ 4 n^{2}}\right)^{(1-2Q) / Q}.
\end{align*}

\end{theorem}
\noindent{\em Proof.}
According to Assumption \ref{assump: lips funon}, we have that
\begin{flalign}\label{Proof Theorem1- 1}
& F_t {  (\boldsymbol{x}_{j,t})} - F_t { (\boldsymbol{x}_{t}^*)} \nonumber \\
& =  F_t {  (\boldsymbol{x}_{j,t})} - F_t { (\boldsymbol{x}_{avg,t})}+ F_t {  (\boldsymbol{x}_{avg,t})} - F_t { (\boldsymbol{x}_{t}^*)} \nonumber \\
&\leq nL_{X}      \|\boldsymbol{x}_{j,t}-\boldsymbol{x}_{avg,t}\|
 +  F_t {  (\boldsymbol{x}_{avg,t})} - F_t{ (\boldsymbol{x}_{t}^*)}\nonumber \\
&\leq   n L_{X} \sum\limits_{i=1}^n   \|\boldsymbol{x}_{i,t}-\boldsymbol{x}_{avg,t}\|  + F_t {  (\boldsymbol{x}_{avg,t})} - F_t{ (\boldsymbol{x}_{t}^*)}.
\end{flalign}

Based on  Algorithm \ref{algorithm 1} and  double stochasticity of $A_{t-1}$, we obtain that ${\boldsymbol{x}}_{avg,t}= \frac{1}{n} \sum_{i=1}^n [\hat{\boldsymbol{x}}_{i,t-1}+ \alpha ( {\boldsymbol{v}}_{i,t-1}-\hat{\boldsymbol{x}}_{i,t-1})]
=\boldsymbol{x}_{avg,t-1} +\alpha ({\boldsymbol{v}}_{avg,t-1}-{\boldsymbol{x}}_{avg,t-1})$.
Thus, by considering Assumption \ref{assump: ball}, we have for any $t\geq2$ that
\begin{flalign}\label{Proof Theorem1- 2-a1}
 \|\boldsymbol{x}_{i,t}-\boldsymbol{x}_{avg,t}\|
&=  \|\hat{\boldsymbol{x}}_{i,t-1}-\boldsymbol{x}_{avg,t-1}+ \alpha (\boldsymbol{v}_{i,t-1}-\hat{\boldsymbol{x}}_{i,t-1})\nonumber \\
&\quad- \alpha (\boldsymbol{v}_{avg,t-1}-{\boldsymbol{x}}_{avg,t-1})\|\nonumber \\
&\leq   \|\hat{\boldsymbol{x}}_{i,t-1}-\boldsymbol{x}_{avg,t-1}\|+ 2\alpha M.
\end{flalign}

With this condition and recalling the regret notion defined in  (\ref{Regret-j}), we obtain from (\ref{Proof Theorem1- 1}) that
\begin{flalign}\label{Proof Theorem1- 1-a2}
&\textbf{Regret}_d^j(T) \nonumber \\
&\leq   n L_{X}
 \sum\limits_{i=1}^n \|\boldsymbol{x}_{i,1}-\boldsymbol{x}_{avg,1}\| +  n L_{X}\sum\limits_{t=1}^{T-1} \sum\limits_{i=1}^n \|\hat{\boldsymbol{x}}_{i,t}-\boldsymbol{x}_{avg,t}\| \nonumber \\
  &\quad +2\alpha Tn^2 L_{X}M
  + \sum\limits_{t=1}^T \left[F_t {  (\boldsymbol{x}_{avg,t})} - F_t{ (\boldsymbol{x}_{t}^*)}\right].
\end{flalign}

Based on this inequality together with the use of Lemma \ref{consistency}, Lemma \ref{grad diffience} and Lemma \ref{lem key regret} in Appendix, we can readily obtain the condition (\ref{Eq Theorem1}) in Theorem \ref{theorem 1}. The proof is completed. \hfill$\square$

Theorem \ref{theorem 1} shows the main results of dynamic regret. It is easy to note that the regret bound of Algorithm \ref{algorithm 1} depends on the choice of $\alpha$. Hence, we have the following corollary by choosing suitable step sizes.
\begin{corollary} \label{corollary 1}
Suppose that the conditions in Theorem \ref{theorem 1} hold.
   Then, if $H_T = o(T), D_T = o(T) $ holds, taking $\alpha = \gamma \sqrt{\frac{H_T +1}{T}}$, we have
\begin{flalign}\label{corollary 1 equation}
 \mathbf{Regret}_d^j(T) \leq \mathcal{O}\left(\sqrt{T(1+H_T)}+D_T\right)
\end{flalign}
where $\gamma$ is a positive adjustment constant guaranteeing $\alpha \leq 1$.
\end{corollary}
\noindent{\em Proof.} According to (\ref{Eq Theorem1}), we obtain
$ \mathbf{Regret}_d^j(T)
\leq \mathcal{O}\left(\alpha T+\frac{1}{\alpha} (1+H_T)+D_T \right).
$
Then, (\ref{corollary 1 equation}) is easily obtained by taking $\alpha T=\frac{1}{\alpha} (1+H_T)$. The proof is complete.
\hfill$\square$

\begin{remark}
In particular, this result shown in Corollary \ref{corollary 1} matches the centralized result in \cite{kalhan2021dynamic} and  is less conservative and tighter than its upper bound $\mathcal{O}(\sqrt{T}(1+H_T+\sqrt{D_T}))$  under the same complexity of linear oracle. Further, the range of $H_T$ in \cite{kalhan2021dynamic} is  improved from $\mathcal{O}(\sqrt{T})$ to $\mathcal{O}(T)$ when a sublinear regret bound is expected, which effectively expands the application field of optimization problems.
\end{remark}
\begin{remark}
Note that $H_T$ as a prior knowledge is difficult to be obtained accurately in practical applications. Two discussions are shown as follows:
\begin{enumerate}
  \item[i)] if an estimated upper bound satisfying $H_T \leq \hat{H}_T=\mathcal{O} (T^{\theta}), 0< \theta < 1 $  can be known in advance, the sublinear dynamic regret is obtained by setting step size $\alpha = \gamma \sqrt{({T^{\theta} +1})/{T}}$.
  \item[ii)] when the loss function $F_t(x)=F(x)$ is time-invariant, $H_T= 0, D_T=0$ hold and the dynamic regret bound $\mathcal{O}(\sqrt{T})$ is established, which matches the results  in distributed off-line scenario, such as \cite{9878290}.
\end{enumerate}
\end{remark}

\section{Simulation}
\begin{figure}[t]
  \centering
\includegraphics[width=8.5cm]{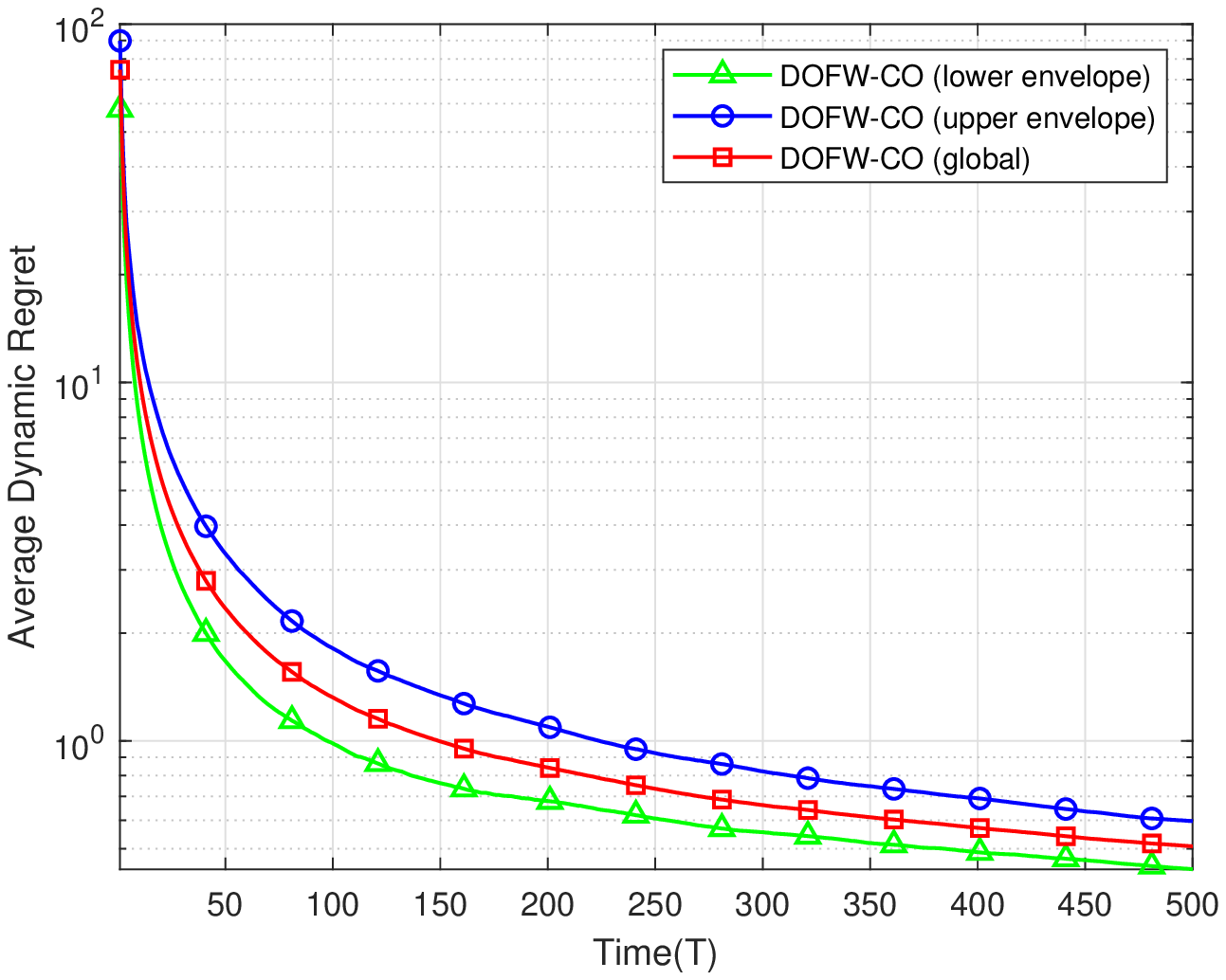}
  \caption{The convergence performance for Algorithm \ref{algorithm 1}.}
  \label{FW_1}
\includegraphics[width=8.5cm]{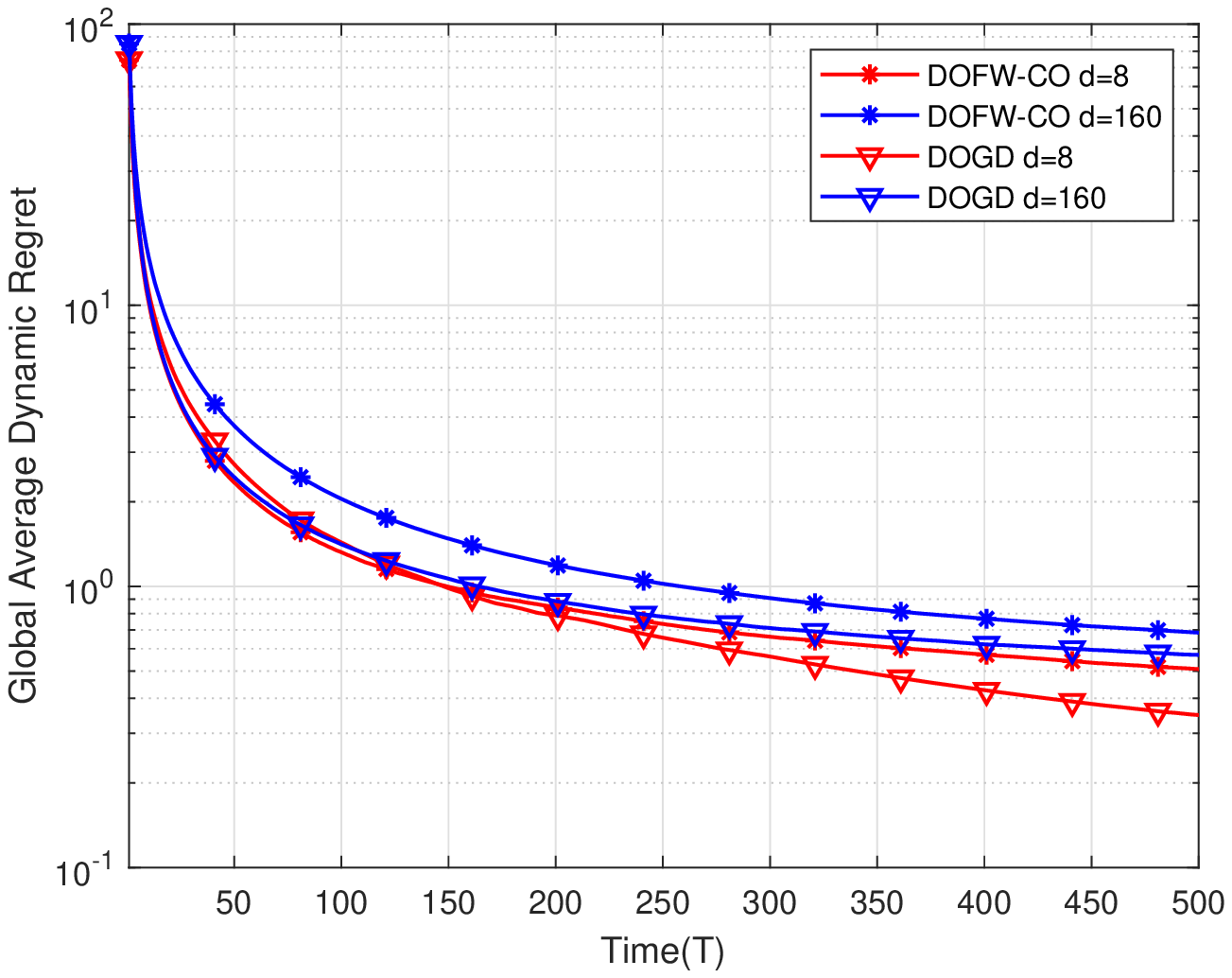}
  \caption{The comparison results of regret between Algorithms \ref{algorithm 1} and DOGD under $d=8$ and $160$.}
  \label{FW_2}
\end{figure}

In this section, several numerical simulations are conducted for ridge regression problem to verify the algorithms we proposed.
The problem of ridge regression is formulated as
\begin{flalign}
&\min\limits_{\boldsymbol{x}\in \boldsymbol{X}} \sum_{t=1}^T \sum_{i=1}^n\left[\frac{1}{2}\left(\boldsymbol{a}_{i, t}^{\top} \bsx -l_{i, t}\right)^2+\lambda_1\|\bsx\|_2^2\right]
\end{flalign}
where $\boldsymbol{X}:=\{\bsx|\  \boldsymbol{1}^T\bsx=1  \}$, $\lambda_1= 5 \times 10^{-6}$ is a regular parameter and the pair $(\boldsymbol{a}_{i, t}, l_{i, t} ) \in \mathbb{R}^d\times\mathbb{R}$ represents the feature and label information only known by agent $i$ at time $t$. The feature vector $\boldsymbol{a}_{i, t}$ is generated randomly and its entries are distributed uniformly from $-5$ to $5$ and the label $l_{i, t}$ satisfies
 $
 l_{i, t}=\boldsymbol{a}_{i, t}^{\top} \bsx_0 + \frac{2\xi_{i,t}}{d\sqrt{t}}
 $
 where $[\bsx_0]_i={1}/d $ and $\xi_{i,t}$ is generated uniformly in the interval $[0, 1]$. We execute  the algorithm over the network of $n=20$ and $d=8$ in this simulation. In the following cases, the global average dynamic regret, the upper envelope and the lower envelope of $ \textbf{Regret}_d^j(T)/T$ are denoted as $\frac{1}{n}\sum_{j=1}^n [\textbf{Regret}_d^j(T)/T]$, $\text{sup}_j\{ \textbf{Regret}_d^j(T)/T\}$ and $\text{inf}_j\{ \textbf{Regret}_d^j(T)/T\}$ to measure the performance of the algorithm, respectively.
 From Fig. \ref{FW_1}, it is clearly seen that the three average dynamic regrets are convergent for Algorithm \ref{algorithm 1} under the condition $\alpha=1 /(4T^{0.4})$, which corresponds to the theoretical result.

\begin{figure}[t]
 \centering
 \includegraphics[width=8.5cm]{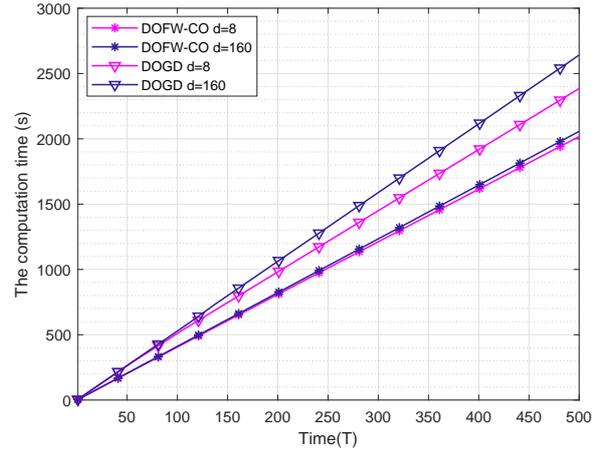}
  \caption{The comparison of the computation time between Algorithms \ref{algorithm 1} and DOGD under $d=8$ and $160$.}
  \label{FW_3}
\end{figure}

  To study the level of computational cost savings of Algorithm \ref{algorithm 1}, we compare the global average dynamic  regret and computational time of Algorithm \ref{algorithm 1} with distributed online gradient descent (DOGD) algorithm \cite{6311406} under two dimensions $d=8$ and $160$, where  the dynamic regret analysis of  Algorithm DOGD can be regarded as a special case satisfying that distance-measuring function  $\mathcal{R}(\bsx)=\| \bsx\|^2$ and the mapping $A=I$ in \cite{shahrampour2017distributed}.  The step sizes $\alpha_1=1 /(4T^{0.4}),\alpha_2=1 /(2T^{0.4})$ and $ \alpha_3=1 /(T^{0.4}),\alpha_4=1 /(4T^{0.4})$ are set for Algorithms \ref{algorithm 1} and DOGD, respectively.
  From Figs. \ref{FW_2} and \ref{FW_3}, we obtain that Algorithm \ref{algorithm 1} is able to achieve similar convergence performance to Algorithm DOGD but using less computation time, which convincingly reflects the advantages of the linear oracle. In particular, when the dimension $d$ is increased from $8$ to $160$, Algorithm \ref{algorithm 1} only has a slightly increase in computation cost compared to the significant increase of  Algorithm DOGD. In term of this point, the higher the dimension of the constrained optimization problem is, the more prominent and important the advantage of Algorithm \ref{algorithm 1} in saving computational costs is.

\section{Conclusions}
For the distributed online convex optimization problem, this paper has developed Algorithm DOFW-CO to reduce the expensive computational cost of the projection step for the high-dimensional and complicated constraint set. In this problem, each agent cooperates to minimize a global cost function through local calculation and information exchange with neighbors over a direct time-varying network.
 We have rigorously established the dynamic regret upper bound of the proposed optimization algorithm as $\mathcal{O}(\sqrt{T(1+H_T)}+D_T)$.
 Finally, the performance of our algorithm has been verified and compared by the simulation of distributed online ridge regression problems. In the future, a promising direction is to explore an improved version of Algorithm DOFW-CO from the perspective of convergence speed or saving computing resources.

\appendix
\section{Some key lemmas}
For the convergence analysis of Algorithm \ref{algorithm 1}, the following lemmas are essential.
Define $\Phi(t, s)=A_{t} A_{t-1} \ldots A_{s}$ as the transition matrix of $A_t$ for all   $t, s \ \text{with} \  t \geq s \geq 1$.

\begin{lemma} \label{consistency}
Let the decision sequence $\{ \bsx_{i,t}\}$ be generated by Algorithm \ref{algorithm 1}. Then, under Assumptions \ref{assump: network} and \ref{assump: ball},  we have for $T\geq 2$ that
\begin{flalign*}
\sumT \sumn \| \hat{\bsx}_{i,t}- \bsx_{avg,t}\| \leq  \frac{n\Gamma}{1-\sigma} \sum\limits_{j=1}^n  \| {\boldsymbol{x}}_{j,1} \|+\alpha T\frac{n^2M\Gamma}{1-\sigma}.
\end{flalign*}
\end{lemma}
\noindent{\em Proof.}
According to Algorithm \ref{algorithm 1}, we get
\begin{flalign} \label{proof_consistency 1}
\hat{\bsx}_{i,t} &=\sum\limits_{j=1}^n [A_t]_{ij}  {\bsx}_{j,t}\\
&= \sum\limits_{j=1}^n [A_t]_{ij}  \hat{\bsx}_{j,t-1} +\alpha \sum\limits_{j=1}^n [A_t]_{ij} (\boldsymbol{v}_{j,t-1}-\hat{\bsx}_{j,t-1} )\nonumber \\
%
&= \sum\limits_{j=1}^n [\Phi(t,1)]_{ij}  {\bsx}_{j,1} +\alpha \sum\limits_{l=1}^{t-1} \sum\limits_{j=1}^n [\Phi(t,l+1)]_{ij} (\boldsymbol{v}_{j,l}-\hat{\bsx}_{j,l} ). \nonumber
\end{flalign}
According to Algorithm \ref{algorithm 1}, the term ${\boldsymbol{x}}_{avg,t}$ can be further simplified as follows:
\begin{flalign} \label{proof_consistency 2}
{\boldsymbol{x}}_{avg,t}&=\frac{1}{n} \sumn \sum_{j=1}^n [A_{t-1}]_{ij} \boldsymbol{x}_{i,t-1}+ \frac{\alpha}{n}\sum_{i=1}^n \left({\boldsymbol{v}}_{i,t-1}-\hat{\boldsymbol{x}}_{j,t-1}\right)\nonumber \\
&={\boldsymbol{x}}_{avg,t-1} + \frac{\alpha}{n}\sum_{i=1}^n \left({\boldsymbol{v}}_{i,t-1}-\hat{\boldsymbol{x}}_{j,t-1}\right)\nonumber \\
&= \frac{1}{n} \sum\limits_{j=1}^n {\boldsymbol{x}}_{j,1} +\frac{\alpha}{n} \sum\limits_{l=1}^{t-1}  \sum\limits_{j=1}^n ({\boldsymbol{v}}_{j,l}-\hat{\boldsymbol{x}}_{j,l})
\end{flalign}
where the third equality combines the double stochasticity of adjacent weight matrix
$A_{t-1}$.

 Combining   (\ref{proof_consistency 1}) and (\ref{proof_consistency 2}), for $t \geq 2$, we achieve
\begin{flalign} \label{proof_consistency 3}
& \| \hat{\bsx}_{i,t}- \bsx_{avg,t}\|\nonumber \\
 &\leq \left\| \sum\limits_{j=1}^n [\Phi(t,1)]_{ij}  {\bsx}_{j,1}- \frac{1}{n} \sum\limits_{j=1}^n {\boldsymbol{x}}_{j,1} \right\|+  \nonumber \\
 &\left\| \alpha \sum\limits_{l=1}^{t-1} \sum\limits_{j=1}^n [\Phi(t,l+1)]_{ij} (\boldsymbol{v}_{j,l}-\hat{\bsx}_{j,l} )- \frac{\alpha}{n} \sum\limits_{l=1}^{t-1}  \sum\limits_{j=1}^n ({\boldsymbol{v}}_{j,l}-\hat{\boldsymbol{x}}_{j,l}) \right\|\nonumber \\
&\leq  \sum\limits_{j=1}^n \left| [\Phi(t,1)]_{ij}  - \frac{1}{n} \right| \| {\boldsymbol{x}}_{j,1} \|\nonumber \\
&+\alpha \sum\limits_{l=1}^{t-1} \sum\limits_{j=1}^n  \left|[\Phi(t,l+1)]_{ij} - \frac{1}{n}\right|   \left\| {\boldsymbol{v}}_{j,l}-\hat{\boldsymbol{x}}_{j,l}\right\|\nonumber \\
&\leq \Gamma \sigma^{t-1}\sum\limits_{j=1}^n  \| {\boldsymbol{x}}_{j,1} \|+\alpha nM\Gamma \sum\limits_{l=1}^{t-1}  \sigma^{t-l-1} \nonumber \\
&\leq \Gamma \sigma^{t-1}\sum\limits_{j=1}^n  \| {\boldsymbol{x}}_{j,1} \|+\alpha \frac{nM\Gamma}{1-\sigma}
\end{flalign}
where the third inequality follows the property of $\Phi(t,s)$ \footnote{\textbf{Lemma}  (\cite{nedic2008distributed})
Let Assumption \ref{assump: network}  hold. Then, for all $i, j \in \mathcal{V}$, we have
$
\left|[\Phi(t, s)]_{i j}-\frac{1}{n}\right| \leq \Gamma \sigma^{(t-s)}
$
where $\sigma=(1-\zeta / 4 n^{2})^{1 /  Q}$ and $\Gamma=(1-\zeta / 4 n^{2})^{(1-2Q) /  Q}$.}
and the fact $\boldsymbol{v}_{j,l}, \hat{\boldsymbol{x}}_{j,l} \in  \boldsymbol{X}$.
Summing from $i=1$ to $n$ and $t=1$ to $T$ on both sides of  (\ref{proof_consistency 3}), we get
\begin{flalign} \label{proof_consistency 4}
 &\sumT \sumn \| \hat{\bsx}_{i,t}- \bsx_{avg,t}\| \nonumber \\
 &\leq \sumn \| \hat{\bsx}_{i,1}- \bsx_{avg,1}\| + n\Gamma \sum_{t=2}^T \sigma^{t-1}\sum\limits_{j=1}^n  \| {\boldsymbol{x}}_{j,1} \|+\alpha T\frac{n^2M\Gamma}{1-\sigma} \nonumber \\
  &\leq    \frac{n\Gamma}{1-\sigma} \sum\limits_{j=1}^n  \| {\boldsymbol{x}}_{j,1} \|+\alpha T\frac{n^2M\Gamma}{1-\sigma}
\end{flalign}
where the second inequality follows the fact $\sum_{i=1}^n \| \hat{\bsx}_{i,1}- \bsx_{avg,1}\| \leq \sum_{i=1}^n \| \sum_{j=1}^n [A_1]_{ij} {\bsx}_{j,1}- \frac{1}{n} \sum_{j=1}^n \bsx_{j,1}\| \leq \sum_{i=1}^n \sum_{j=1}^n| [A_1]_{ij}- \frac{1}{n}| \|  \bsx_{j,1}\| \leq n\Gamma  \sum_{j=1}^n  \|  \bsx_{j,1}\|$.
The proof is complete.
\hfill$\square$
\begin{lemma} \label{Sigma Tracking diffience}
 Let the decision sequence $\{ \bsx_{i,t}\}$ be generated by Algorithm \ref{algorithm 1}. Then, under Assumptions \ref{assump: ball} and \ref{assump: lips grad},  we have  for any $T\geq 2$ that
\begin{flalign}
\textbf{(a)} \ &\sumn \overline{\nabla} f_{i,t} = \sumn \nabla f_{i,t} (\hat{\boldsymbol{x}}_{i,t}),\ \forall t \geq 1. \label{a-lem3}\\
\textbf{(b)} \  &\sum\limits_{t=2}^{T} \sumn \left\|  \boldsymbol{\delta}_{i,t}  \right\| \leq2G_{X}\sum\limits_{t=2}^{T}\sumn  \| \hat{\boldsymbol{x}}_{i,t-1} -{\boldsymbol{x}}_{avg,t-1} \|     \nonumber \\
 &\quad \quad \quad \quad \quad \quad \quad +  n D_T   +n M G_{X} \alpha T. \label{b-lem3}
\end{flalign}
\end{lemma}
\noindent{\em Proof.}
\textbf{(a)} We utilize the mathematical induction method to prove this part. Note that $\overline{\nabla} f_{i,1} = \nabla f_{i,1} (\hat{\boldsymbol{x}}_{i,1})$ according to Algorithm \ref{algorithm 1}. Thus,   the equality in (\ref{a-lem3})  holds when $t=1$. Now  we assume that $\sum_{i=1}^n \overline{\nabla} f_{i,t} = \sum_{i=1}^n \nabla f_{i,t} (\hat{\boldsymbol{x}}_{i,t})$ holds at some $t$, and we are going to show this equality also holds at $t+1$. Actually,
%
%
\begin{flalign}
& \sumn \overline{\nabla} f_{i,t+1}\nonumber \\
& = \sumn \widehat{\nabla} f_{i,t} + \sumn  \nabla f_{i,t+1}  (\hat{\boldsymbol{x}}_{i,t+1}) - \sumn  \nabla f_{i,t}  (\hat{\boldsymbol{x}}_{i,t}) \nonumber \\
&= \sumn \sum_{j=1}^n [A_t]_{ij} \overline{\nabla} f_{j,t} + \sumn  \nabla f_{i,t+1}  (\hat{\boldsymbol{x}}_{i,t+1}) - \sum_{i=1}^n \overline{\nabla} f_{i,t}  \nonumber \\
&=\sumn  \nabla f_{i,t+1}  (\hat{\boldsymbol{x}}_{i,t+1})
\end{flalign}
where the last equality follows from the double stochasticity of $A_{t}$.

\textbf{(b)} By using (\ref{average-delta}), we obtain for any $t\geq2$ that
\begin{flalign} \label{proof Sigma Tracking 1}
 \sumn \left\|  \boldsymbol{\delta}_{i,t}  \right\| & \leq \sumn \left\| \nabla f_{i,t}  (\hat{\bsx}_{i,t}) - \nabla f_{i,t-1}  (\hat{\boldsymbol{x}}_{i,t})\right\| \nonumber \\
 &\quad +\sumn \left\| \nabla f_{i,t-1}  (\hat{\boldsymbol{x}}_{i,t}) - \nabla f_{i,t-1}  (\hat{\boldsymbol{x}}_{i,t-1})\right\| \nonumber \\
 &\leq n g_{t-1,sup} +\sumn \left(  G_{X} \| \hat{\boldsymbol{x}}_{i,t} -\hat{\boldsymbol{x}}_{i,t-1} \|  \right)\nonumber \\
 &\leq n g_{t-1,sup} +G_{X}\sumn (  \| \hat{\boldsymbol{x}}_{i,t} -{\boldsymbol{x}}_{avg,t-1} \|  \nonumber \\
 &\quad +\| \hat{\boldsymbol{x}}_{i,t-1}-{\boldsymbol{x}}_{avt,t-1} \| )\nonumber \\
 &\leq n g_{t-1,sup} + 2G_{X}\sumn \| \hat{\boldsymbol{x}}_{i,t-1} -{\boldsymbol{x}}_{avg,t-1} \|   \nonumber \\
 & \quad +n M G_{X} \alpha
\end{flalign}
where the last inequality is obtained based on the fact:
$
\sum_{i=1}^n \| \hat{\boldsymbol{x}}_{i,t} -{\boldsymbol{x}}_{avg,t-1} \|
\leq \sum_{i=1}^n \sum_{j=1}^n [A_t]_{ij}\| {\boldsymbol{x}}_{j,t} -{\boldsymbol{x}}_{avg,t-1} \|
\leq \sum_{i=1}^n \| \hat{\boldsymbol{x}}_{i,t-1} -{\boldsymbol{x}}_{avg,t-1} +\alpha (\boldsymbol{v}_{i,t-1}-\hat{\boldsymbol{x}}_{i,t-1})\|
\leq \sum_{i=1}^n \| \hat{\boldsymbol{x}}_{i,t-1} -{\boldsymbol{x}}_{avg,t-1} \| + \alpha n M.
$

Then, by summing the both sides of  (\ref{proof Sigma Tracking 1}) from  $t=2$ to $T$, we can readily obtain the inequality in (\ref{b-lem3}). The proof is complete.
 \hfill$\square$

\begin{lemma} \label{grad diffience}
Let the sequence $\{ \widehat{\nabla} f_{i,t}, {\nabla} f_{i,t} (\hat{\boldsymbol{x}}_{i,t}) \}$ be generated by Algorithm \ref{algorithm 1}. Then, under Assumptions \ref{assump: network} and \ref{assump: lips grad}, we have for any $T\geq 2$ that
\begin{flalign} \label{condition-lem4}
& \sumT \sumn \left\|  \widehat{\nabla} f_{i,t} - \frac{1}{n} \sum_{j=1}^n   {\nabla} f_{j,t} (\hat{\boldsymbol{x}}_{j,t})\right\| \nonumber \\
& \leq  \frac{ n \Gamma}{1-\sigma} \sum\limits_{j=1}^n\left\|   {\nabla} f_{j,1} (\hat{\boldsymbol{x}}_{j,1}) \right\|  + \frac{ n \Gamma}{1-\sigma} \sum\limits_{t=2}^{T} \sumn \left\|  \boldsymbol{\delta}_{i,t}  \right\|.
\end{flalign}
\end{lemma}
\noindent{\em Proof.} Similar to the proof of Lemma \ref{consistency} and combining  Algorithm \ref{algorithm 1}, (\ref{a-lem3}),   for any  $t\geq2$ it can be verified that
\begin{flalign} \label{proof_grad diff 1}
&\widehat{\nabla} f_{i,t} 
%
%
= \sum\limits_{j=1}^n [\Phi(t,2)]_{ij}\widehat{\nabla} f_{j,1} + \sum\limits_{l=2}^t \sum\limits_{j=1}^n [\Phi(t,l)]_{ij} \boldsymbol{\delta}_{j,l},\\
%
&\frac{1}{n} \sumn   {\nabla} f_{i,t} (\hat{\boldsymbol{x}}_{i,t})= \frac{1}{n} \sumn   \overline{\nabla} f_{i,1} + \frac{1}{n} \sum\limits_{l=2}^t\sumn  \boldsymbol{\delta}_{i,l}. \label{proof_grad diff 2}
\end{flalign}

Similar to (\ref{proof_consistency 3}), combining the fact $\widehat{\nabla} f_{i,1}=\sum_{j=1}^n [A_1]_{ij}  \overline{\nabla} f_{j,1},  \overline{\nabla} f_{i,1}= {\nabla} f_{i,1} (\hat{\boldsymbol{x}}_{i,1})$, it follows from (\ref{proof_grad diff 1}) and (\ref{proof_grad diff 2}) that $\|  \widehat{\nabla} f_{i,t} - \frac{1}{n} \sum_{j=1}^n  {\nabla} f_{j,t} (\hat{\boldsymbol{x}}_{j,t})\| \leq\sum_{j=1}^n \Gamma \sigma^{t-1}  \|   {\nabla} f_{j,1} (\hat{\boldsymbol{x}}_{j,1}) \| +\sum_{l=2}^t \sum_{j=1}^n\Gamma \sigma^{t-l} \|  \boldsymbol{\delta}_{j,l} \|$.
This implies that
\begin{flalign} \label{proof_grad diff 4}
 &\sumT \sumn \left\|  \widehat{\nabla} f_{i,t} - \frac{1}{n} \sum_{j=1}^n   {\nabla} f_{j,t} (\hat{\boldsymbol{x}}_{j,t})\right\|\nonumber \\
 &\leq\sumn \left\|  \widehat{\nabla} f_{i,1} -\frac{1}{n} \sum_{j=1}^n    {\nabla} f_{j,1} (\hat{\boldsymbol{x}}_{j,1})\right\| \nonumber \\
 &+\sum\limits_{t=2}^T \sum\limits_{j=1}^n n \Gamma \sigma^{t-1}  \left\|   {\nabla} f_{j,1} (\hat{\boldsymbol{x}}_{j,1}) \right\| +\sum\limits_{t=2}^T \sum\limits_{l=2}^t \sum\limits_{j=1}^n n \Gamma \sigma^{t-l} \left\|  \boldsymbol{\delta}_{j,l}  \right\| \nonumber \\
 &\leq\sumn \sum\limits_{j=1}^n \left | [A_1]_{ij} -  \frac{1}{n} \right| \left\|  {\nabla} f_{j,1} (\hat{\boldsymbol{x}}_{j,1}) \right\| \nonumber \\
 &+\frac{\sigma n \Gamma}{1-\sigma} \sum\limits_{j=1}^n \left\|   {\nabla} f_{j,1} (\hat{\boldsymbol{x}}_{j,1}) \right\| +n\Gamma \left(\sum\limits_{l=0}^{T-2}
 \sigma^l \right)\left(\sum\limits_{t=2}^{T} \sumn \left\|  \boldsymbol{\delta}_{i,t}  \right\| \right)\nonumber \\
 &\leq \frac{n \Gamma}{1-\sigma} \sum\limits_{j=1}^n \left\|   {\nabla} f_{j,1} (\hat{\boldsymbol{x}}_{j,1}) \right\|+\frac{ n \Gamma}{1-\sigma} \sum\limits_{t=2}^{T} \sumn \left\|  \boldsymbol{\delta}_{i,t}  \right\|
\end{flalign}
Substituting the above inequalities into   (\ref{proof_grad diff 4}),  we can readily obtain (\ref{condition-lem4}). The proof is complete.
\hfill$\square$
\begin{lemma} \label{lem key regret}
Let the decision sequence $\{ \bsx_{i,t}\}$ be generated by Algorithm \ref{algorithm 1}. Then, under Assumptions \ref{assump: ball}, \ref{assump: lips funon} and \ref{assump: lips grad},  we have  for any $T\geq 2$ that
\begin{flalign} \label{eq lem key regret}
  &\sum\limits_{t=1}^{T}  F_{t} (\bsx_{avg,t})- \sum\limits_{t=1}^{T} F_{t}(\bsx_t^*) \nonumber \\
  &\leq \frac{2n}{\alpha} H_T + \frac{nL_X M}{\alpha}  + { 2M}\sum\limits_{t=1}^{T-1}  \sumn \left\| \sum_{j=1}^n   {\nabla} f_{j,t} (\hat{\boldsymbol{x}}_{j,t})- \widehat{\nabla} f_{i,t} \right\|\nonumber \\
 &\quad + {2 M G_X} \sum_{t=1}^{T-1} \sumn \| \hat{\bsx}_{i,t}-\bsx_{avg,t}\| + \frac{ nG_X M^2}{2} \alpha T+   n M  D_T.
\end{flalign}

\end{lemma}

\noindent{\em Proof.}
By using the smooth property in Assumption \ref{assump: lips grad}, we have
\begin{flalign} \label{proof lem key 1}
&F_{t+1} (\boldsymbol{x}_{avg,t+1}) -F_{t+1} (\boldsymbol{x}_{avg,t})  \\
&\leq  \left< \nabla F_{t+1} (\boldsymbol{x}_{avg,t}), \boldsymbol{x}_{avg,t+1}- \boldsymbol{x}_{avg,t} \right> \nonumber \\
 &\quad+ \frac{nG_X}{2} \| \boldsymbol{x}_{avg,t+1}- \boldsymbol{x}_{avg,t} \|^2 \nonumber \\
&\leq \alpha \sumn \left< \frac{1}{n}\nabla F_{t+1} (\boldsymbol{x}_{avg,t}), \boldsymbol{v}_{i,t}- \boldsymbol{x}_{avg,t} \right> + \frac{ nG_X M^2\alpha^2}{2}. \nonumber
\end{flalign}
It can be further verified that
\begin{flalign} \label{proof lem key 2}
& \left< \frac{1}{n}\nabla F_{t+1} (\boldsymbol{x}_{avg,t}), \boldsymbol{v}_{i,t}- \boldsymbol{x}_{avg,t} \right> \nonumber \\
&\leq  \left< \frac{1}{n}\nabla F_{t+1} (\boldsymbol{x}_{avg,t})- \widehat{\nabla} f_{i,t}, \boldsymbol{v}_{i,t}- \boldsymbol{x}_{avg,t} \right> \nonumber \\
&\quad + \left< \widehat{\nabla} f_{i,t}, \boldsymbol{x}_t^*- \boldsymbol{x}_{avg,t} \right> \nonumber \\
&\leq \frac{ M}{n} \left\|  \nabla F_{t+1} (\bsx_{avg,t}) - \nabla F_{t} (\bsx_{avg,t}) \right\|+{ 2M} \left\| \frac{1}{n}\nabla F_{t} (\bsx_{avg,t}) \right.\nonumber \\
&\quad \left. -\frac{1}{n} \sum_{j=1}^n   {\nabla} f_{j,t} (\hat{\boldsymbol{x}}_{j,t}) +\frac{1}{n} \sum_{j=1}^n   {\nabla} f_{j,t} (\hat{\boldsymbol{x}}_{j,t})- \widehat{\nabla} f_{i,t} \right\| \nonumber \\
   &\quad + \frac{1}{n} \left< \nabla F_{t} (\bsx_{avg,t}), \bsx_t^*- \bsx_{avg,t} \right>
   \nonumber \\
&\leq  { M} g_{t,sup} +{ 2M}  \left\| \frac{1}{n} \sum_{j=1}^n   {\nabla} f_{j,t} (\hat{\boldsymbol{x}}_{j,t})- \widehat{\nabla} f_{i,t} \right\|  \\
&\quad+\frac{2M G_X}{n} \sum _{j=1}^n  \| \hat{\bsx}_{j,t}-\bsx_{avg,t}\|
   +  \frac{1}{n}\left[ F_{t}(\bsx_t^*) - F_{t} (\bsx_{avg,t}) \right] \nonumber
\end{flalign}
where the first inequality is obtained by utilizing the following  optimality condition:
\begin{flalign}
\left<\widehat{\nabla} f_{i,t}, \boldsymbol{x}_t^*\right>  \geq \min\limits_{\bsx  \in \boldsymbol{X}} \left< \widehat{\nabla} f_{i,t}, \boldsymbol{x} \right>  =  \left< \widehat{\nabla} f_{i,t}, \boldsymbol{v}_{i,t} \right>
\end{flalign}
and the last inequality is derived  based on  the convexity condition of $F_t( \bsx)$ together with   Assumption \ref{assump: lips grad}. Then, it follows from (\ref{proof lem key 1}) and  (\ref{proof lem key 2}) that
\begin{flalign} \label{proof lem key 3}
&F_{t+1} (\boldsymbol{x}_{avg,t+1})-F_{t+1} (\boldsymbol{x}_{avg,t}) \nonumber \\
&\leq  {  \alpha n M} g_{t,sup} + \alpha \left[ F_{t}(\bsx_t^*) - F_{t} (\bsx_{avg,t}) \right] +\Psi_t.
\end{flalign}
where $\Psi_t={ 2 \alpha M}  \sum_{i=1}^n \left\| \frac{1}{n} \sum_{j=1}^n   {\nabla} f_{j,t} (\hat{\boldsymbol{x}}_{j,t})- \widehat{\nabla} f_{i,t} \right\|+ {2\alpha M G_X} \sum_{j=1}^n \| \hat{\bsx}_{j,t}-\bsx_{avg,t}\|+ \frac{ nG_X M^2\alpha^2}{2}$.
Based on this inequality, we can further obtain that
\begin{flalign} \label{proof lem key 4}
&F_{t+1} (\boldsymbol{x}_{avg,t+1})-F_{t+1}(\boldsymbol{x}_{t+1}^*)\nonumber \\
&\leq F_{t+1} (\boldsymbol{x}_{avg,t}) -F_{t+1}(\boldsymbol{x}_{t+1}^*)
 +{  \alpha n M} g_{t,sup} +\Psi_t \nonumber \\
&\quad + \alpha \left[ F_{t}(\bsx_t^*) - F_{t} (\bsx_{avg,t}) \right] \nonumber \\
&\leq n f_{t,sup}+F_{t} (\boldsymbol{x}_{avg,t})- F_{t}(\boldsymbol{x}_{t}^*) +F_{t}(\boldsymbol{x}_{t}^*) -F_{t+1}(\boldsymbol{x}_{t+1}^*)\nonumber \\
&\quad  +{  \alpha n M} g_{t,sup}+ \alpha \left[ F_{t}(\bsx_t^*) - F_{t} (\bsx_{avg,t}) \right] + \Psi_t
\end{flalign}
where the last inequality is established by using the following fact:
\begin{flalign*}
&F_{t+1} (\boldsymbol{x}_{avg,t}) -F_{t+1}(\boldsymbol{x}_{t+1}^*) \nonumber \\
&= F_{t+1} (\boldsymbol{x}_{avg,t}) - F_{t} (\boldsymbol{x}_{avg,t}) +F_{t} (\boldsymbol{x}_{avg,t})- F_{t}(\boldsymbol{x}_{t}^*)\nonumber \\
&\quad +F_{t}(\boldsymbol{x}_{t}^*) -F_{t+1}(\boldsymbol{x}_{t+1}^*) \nonumber \\
&\leq n f_{t,sup}+F_{t} (\boldsymbol{x}_{avg,t})- F_{t}(\boldsymbol{x}_{t}^*) +F_{t}(\boldsymbol{x}_{t}^*) -F_{t+1}(\boldsymbol{x}_{t+1}^*).
\end{flalign*}

By summing the both sides of (\ref{proof lem key 4}), we get
\begin{flalign} \label{proof lem key 6}
&\sum\limits_{t=1}^{T-1} \left\{ F_{t+1} (\boldsymbol{x}_{avg,t+1})-F_{t+1}(\boldsymbol{x}_{t+1}^*) - [F_{t} (\boldsymbol{x}_{avg,t})- F_{t}(\boldsymbol{x}_{t}^*)] \right\} \nonumber \\
&\quad + \alpha \sum\limits_{t=1}^{T-1} \left[ F_{t} (\bsx_{avg,t})-F_{t}(\bsx_t^*) \right] \nonumber \\
&\leq n H_T +F_{1}(\boldsymbol{x}_{1}^*) -F_{T}(\boldsymbol{x}_{T}^*)  +{  \alpha n M} D_T + \sum\limits_{t=1}^{T-1} \Psi_t
\end{flalign}

By recalling the facts that $0< \alpha \leq1$ and $F_{T} (\boldsymbol{x}_{avg,T})-F_{T}(\boldsymbol{x}_{T}^*)\geq 0$, we have  that
\begin{flalign} \label{proof lem key 6-a2}
&\alpha \sum\limits_{t=1}^{T} \left[ F_{t} (\bsx_{avg,t})-F_{t}(\bsx_t^*) \right] \\
&\leq [F_{T} (\boldsymbol{x}_{avg,T})-F_{T}(\boldsymbol{x}_{T}^*)] + \alpha \sum\limits_{t=1}^{T-1} \left[ F_{t} (\bsx_{avg,t})-F_{t}(\bsx_t^*) \right] \nonumber \\
&\leq n H_T + [F_{1} (\boldsymbol{x}_{avg,1})-F_{T}(\boldsymbol{x}_{T}^*)] +{  \alpha n M} D_T + \sum\limits_{t=1}^{T-1} \Psi_t. \nonumber
\end{flalign}

In addition, it can be obtained under Assumption \ref{assump: ball} and Assumption \ref{assump: lips funon} that
$
F_{1} (\boldsymbol{x}_{avg,1}) -F_{T}(\boldsymbol{x}_{T}^*)
=F_{1} (\boldsymbol{x}_{avg,1})-F_{1}(\boldsymbol{x}_{T}^*)+\sum_{t=1}^{T-1} [F_{t}(\boldsymbol{x}_{T}^*)-F_{t+1}(\boldsymbol{x}_{T}^*)]
\leq nL_X \|\boldsymbol{x}_{avg,1}- \boldsymbol{x}_{T}^*\|+\sum_{t=1}^{T-1} \sum_{i=1}^n \left|f_{i,t}(\boldsymbol{x}_{T}^*)-f_{i,t+1}(\boldsymbol{x}_{T}^*)\right|
\leq nL_X M +n H_T$.
%
%
%

This, together with (\ref{proof lem key 6-a2}), yields that
\begin{flalign} \label{proof lem key 6-a6}
&\alpha \sum\limits_{t=1}^{T} \left[ F_{t} (\bsx_{avg,t})-F_{t}(\bsx_t^*) \right] \nonumber \\
&\leq 2n H_T +nL_X M +  \alpha n M  D_T + \sum\limits_{t=1}^{T-1} \Psi_t.
\end{flalign}
This further implies the result in (\ref{eq lem key regret}) through recalling the definition of $\Psi_t$.  The proof is complete.
\hfill$\square$

\bibliographystyle{model5-names}
\bibliography{AutoV3_DR_DOFW_CO}

\end{document}